# High Dimensional Probability

## Proceedings of the Fourth International Conference

Evarist Giné, Vladimir Koltchinskii, Wenbo Li, Joel Zinn, Editors

Institute of Mathematical Statistics
*Lecture Notes–Monograph Series*

Series Editor:
R. A. Vitale

The production of the *Institute of Mathematical Statistics Lecture Notes–Monograph Series* is managed by the IMS Office: Jiayang Sun, Treasurer and Elyse Gustafson, Executive Director.

Library of Congress Control Number: 2001012345

International Standard Book Number 978-0-940600-67-6,

0-940600-67-6

# Contents

## APPLICATIONS OF EMPIRICAL PROCESSES

# Preface

About forty years ago it was realized by several researchers that the essential features of certain objects of Probability theory, notably Gaussian processes and limit theorems, may be better understood if they are considered in settings that do not impose structures extraneous to the problems at hand. For instance, in the case of sample continuity and boundedness of Gaussian processes, the essential feature is the metric or pseudometric structure induced on the index set by the covariance structure of the process, regardless of what the index set may be. This point of view ultimately led to the Fernique-Talagrand majorizing measure characterization of sample boundedness and continuity of Gaussian processes, thus solving an important problem posed by Kolmogorov. Similarly, separable Banach spaces provided a minimal setting for the law of large numbers, the central limit theorem and the law of the iterated logarithm, and this led to the elucidation of the minimal (necessary and/or sufficient) geometric properties of the space under which different forms of these theorems hold. However, in light of renewed interest in Empirical processes, a subject that has considerably influenced modern Statistics, one had to deal with a non-separable Banach space, namely $\mathcal{L}_\infty$. With separability discarded, the techniques developed for Gaussian processes and for limit theorems and inequalities in separable Banach spaces, together with combinatorial techniques, led to powerful inequalities and limit theorems for sums of independent bounded processes over general index sets, or, in other words, for general empirical processes.

This research led to the introduction or to the re-evaluation of many new tools, including randomization, decoupling, chaining, concentration of measure and exponential inequalities, series representations, that are useful in other areas, among them, asymptotic geometric analysis, Banach spaces, convex geometry, nonparametric statistics, computer science (e.g. learning theory).

The term High Dimensional Probability, and Probability in Banach spaces before, refers to research in probability and statistics that emanated from the problems mentioned above and the developments that resulted from such studies.

A large portion of the material presented here is centered on these topics. For example, under limit theorems one has represented both the theoretical side as well as applications to Statistics; research on dependent as well as independent random variables; Lévy processes as well as Gaussian processes; U and V-processes as well as standard empirical processes. Examples of tools to handle problems on such topics include concentration inequalities and stochastic inequalities for martingales and other processes. The applications include classical statistical problems and newer areas such as Statistical Learning theory.

Many of the papers included in this volume were presented at the IVth International Conference on High Dimensional Probability held at St. John's College, Santa Fe, New Mexico, on June 20-24, 2005, and all of them are based on topics covered at this conference. This conference was the fourteenth in a series that began with the Colloque International sur les Processus Gaussiens et les Distributions Aléatoires, held in Strasbourg in 1973, continued with nine conferences on Probability in Banach Spaces, and four with the title of High Dimensional Probability. The book *Probability in Banach Spaces* by M. Ledoux and M. Talagrand, Springer-Verlag 1991, and the Preface to the volume *High Dimensional Probability III*, Birkhäuser, 2003, contain information on these precursor conferences. More historical information can be found online at http://www.math.udel.edu/˜wli/hdp/index.html. This last reference also includes a list of titles of talks and participants of this meeting.

The participants to this conference are grateful for the support of the National Science Foundation, the National Security Agency and the University of New Mexico.

July, 2006

Evarist Giné
Vladimir Koltchinskii
Wenbo Li
Joel Zinn

# Contributors to this volume

Arcones, M. A., *Binghamton University*

Berthet, P., *Université Rennes 1*
Borisov, I. S., *Sobolev Institute of Mathematics*
Bystrov, A. A., *Sobolev Institute of Mathematics*

Deheuvels, P., *L.S.T.A., Université Pierre et Marie Curie (Paris 6)*
Dobrić, V., *Lehigh University*
Dudley, R. M., *Massachusetts Institute of Technology*

Eggermont, P. P. B., *University of Delaware*

Figueroa-López, J. E., *Purdue University*

Giné, E., *University of Connecticut*

Houdré, C., *Georgia Institute of Technology*
Hush, D., *Los Alamos National Laboratory*

Koltchinskii, V., *Georgia Institute of Technology*

LaRiccia, V. N., *University of Delaware*

Mason, D. M., *University of Delaware*
Mendelson, S., *ANU & Technion I.I.T*

Nickl, R., *Department of Statistics, University of Vienna*

Ojeda, F. M., *Universidad Simón Bolívar and Lehigh University*

Peligrad, M., *University of Cincinnati*
Pinelis, I., *Michigan Technological University*

Scovel, C., *Los Alamos National Laboratory*
Steinberg, S., *University of New Mexico*
Steinwart, I., *Los Alamos National Laboratory*

Umarov, S., *University of New Mexico*
Utev, S., *University of Nottingham*

Vershynin, R., *University of California, Davis*

Wu, D., *Department of Statistics*
Wu, W. B., *University of Chicago*

Xiao, Y., *Probability, Michigan State University*

Zinn, J., *Texas A&M University*